\documentclass[12pt]{amsart}

\usepackage{amsmath, amssymb, amsfonts}

\newcommand {\ev} {{\bar0}}
\newcommand {\od} {{\bar1}}
\newcommand {\eps} {\varepsilon}

\newcommand {\cal} {\mathcal}

\newcommand {\cD}     {{\cal D}}

\newcommand {\Cee}    {{\mathbb  C}}

\newcommand {\Ree}    {{\mathbb  R}}
\newcommand {\Zee}    {{\mathbb  Z}}
\newcommand {\Kee}    {{\mathbb  K}}

\newcommand {\fa}     {{\mathfrak{a}}}

\newcommand {\fas}    {{\mathfrak{as}}}

\newcommand {\fc}    {{\mathfrak{c}}}
\newcommand {\fp}    {{\mathfrak{p}}}
\newcommand {\fs}    {{\mathfrak{s}}}

\newcommand {\fder}   {{\mathfrak{der}}}   %

\newcommand {\fd}     {{\mathfrak{d}}}
\newcommand {\fe}     {{\mathfrak{e}}}

\newcommand {\ff}     {{\mathfrak{f}}}

\newcommand {\fg}     {{\mathfrak{g}}}    %
\newcommand {\fgl}    {{\mathfrak{gl}}}  %
\newcommand {\fh}     {{\mathfrak{h}}}

\newcommand {\fii}    {{\mathfrak{i}}}    %
\newcommand {\fk}     {{\mathfrak{k}}}

\newcommand {\fn}     {{\mathfrak{n}}}

\newcommand {\fo}     {{\mathfrak{o}}}

\newcommand {\fsl}    {{\mathfrak{sl}}}

\newcommand {\fsp}    {{\mathfrak{sp}}}

\newcommand {\fsvect} {{\mathfrak{svect}}}
\newcommand {\fvect}  {{\mathfrak{vect}}}   %
\newcommand {\fv}   {{\mathfrak{v}}}

\newcommand {\subplus}{\mathop{{\subset}\llap{\raise
0.5pt\hbox{\normalfont\small+}\hskip 0.5pt}}}

\newcommand{\rmname}[1]
  {\expandafter\newcommand \csname #1\endcsname {{\operatorname{#1}}}}

\newcommand{\rmnameii}[2]
  {\expandafter\newcommand \csname #1\endcsname {{\operatorname{#2}}}}

\rmname{ad} \rmname{const} \rmname{End} \rmname{id} \rmname{Span}
\rmname{Pty}

\rmnameii{Div}{div}
\rmnameii {Char} {char}
\rmnameii {IM} {Im}
\rmnameii {vvol} {vol}

\newcommand {\tto} {\longrightarrow}

\newtheorem{Theorem}{Theorem}[section]
\newtheorem{Lemma}[Theorem]{Lemma}

\theoremstyle{remark}
\newtheorem{Remark}[Theorem]{Remark}

\begin{document}

\title[Lie algebras of vector fields]{How to realize Lie algebras by vector fields}

\author{Irina Shchepochkina}

\address{Independent University of Moscow,
Bolshoj Vlasievsky per, dom 11, RU-121 002 Moscow, Russia;
irina@mccme.ru}



\begin{abstract} An algorithm for embedding finite dimensional Lie algebras
into Lie algebras of vector fields (and Lie superalgebras into Lie
superalgebras of vector fields) is offered in a way applicable
over ground fields of any characteristic. The algorithm is
illustrated by reproducing Cartan's interpretations of the Lie
algebra of G(2) as the Lie algebra that preserves certain
non-integrable distributions. Similar algorithm and interpretation
are applicable to other exceptional simple Lie algebras, as well
as to all non-exceptional simple ones and many non-simple ones,
and to many Lie superalgebras.
\end{abstract}

\thanks{Financial support of RFBR grant 05-01-00001 is gratefully acknowledged.}

\maketitle

\section{Introduction}

\begin{quote}\rightline{In memory of Felix Aleksandrovich Berezin}\end{quote}

Here I offer an algorithm which explicitly describes how to embed
any $\Zee$-graded Lie algebra (or Lie superalgebra)
$\fn:=\mathop{\oplus}\limits_{k\geq-d}\fn_k$ such that
\begin{equation}
\label{1} \fn_{-1}\;\text{ generates
$\fn_-:=\mathop{\oplus}\limits_{k<0}\fn_k$ and $\dim
\fn_-<\infty$}
\end{equation}
into a Lie algebra (resp., Lie superalgebra) of polynomial vector
fields over $\Ree$ or $\Cee$ or over a field $\Kee$ of
characteristic $p>0$.\footnote{Although $p$ denotes the
characteristic of the ground field, parity, and is used as an
index, the context is always clear. }

For almost a decade, whenever asked, I described the  algorithm I
propose here but was reluctant to publish it as a research paper:
the algorithm is straightforward and was, actually, used more than
a century ago by Cartan \cite{C}, and recently by Yamaguchi
\cite{Y}. For the same reason, the draft of \cite{ShE} with some
examples of embeddings based on this algorithm was also being put
aside and will appear as a sequel to this paper; in the meantime
\cite{La1} and \cite{La2} appeared with some more
examples\footnote{I was unable to follow the details of both
\cite{La1} and \cite{La2} and hopefully this paper will help to
elucidate important realizations of \cite{La1} and \cite{La2}.}.

Grozman and Leites convinced me, however, that the algorithm, and
its usefulness, were never expressed explicitly. Most
convincingly, they used the algorithm not only for interpreting
known, but mysterious, simple Lie algebras, and Lie superalgebras,
especially in characteristic $p>0$, but in order to get new
examples in the absence of classification (\cite{GL}). So here it
is.  Grozman already implemented it in his {\bf SuperLie} package
\cite{Gr}.

Having started to write, I added something new as compared with
\cite{C}: a description by means of differential equations of {\it
partial} prolongs --- subalgebras of the Lie algebras of
polynomial vector fields embedded \lq\lq projective-like\rq\rq.
Such description is particularly important if $p>0$, and for some
Lie superalgebras.

At the last moment, I learned that, for $p>0$, Fei and Shen
\cite{FSh} proved existence of embeddings I consider and
illustrated it with a description of the simple Lie algebras of
contact vector fields for $p=2$. They also formulate questions
this paper answers.

For reviews of related to our result realizations of Lie
(super)algebras by differential operators (not necessarily first
order homogeneous ones), see \cite{BGLS, M} and \cite{Sh2}.

{\bf Problem formulation,  facts known, and our reasons}. Let
$\fn:=\mathop{\oplus}\limits_{k=-d}^{-1}\fn_k$, be an
$n$-dimensional $\Zee$-graded Lie algebra of depth $d>1$
satisfying (\ref{1}). Let $f:\fn\tto\fvect(n)=\fder\Kee[x_1, \dots
, x_n]$ be an embedding. The image $f(\fn)$ is a subspace in the
space of vector fields; every vector field can be evaluated at any
point; let $f(\fn)(0)$ be the span of these evaluations at $0$.

{\bf Problem 1}. {\it Embed $\fn$ into $\fvect(n)$ so that the
$\dim f(\fn)(0)=n$. }

{\bf Comment}. Roughly speaking, we wish the image of $\fn$ be
spanned by all partial derivatives modulo vector fields that
vanish at the origin.

 \vskip 0.2 cm

Such an embedding determines a non-standard\footnote{The grading
$\deg x_i=1$ for all $i$ associated with the $(x)$-adic filtration
is said to be {\it standard}; any other grading is {\it
non-standard}.} grading of depth $d$ on $\fvect(n)$. We will
denote $\fvect(n)$ considered with this non-standard grading by
$\fv= \mathop{\oplus}\limits_{k=-d}^{\infty}\fv_k$. Let $\fg_-$ be
the image of $\fn$ in $\fv$, i.e., $\fg_-\subset
\fv_-:=\mathop{\oplus}\limits_{k<0}\fv_k$.

\vskip 0.2 cm

{\bf Problem 2}.  {\it Compute the {\bf complete algebraic
prolong} of $\fg_-$, i.e., the maximal subalgebra
$(\fg_-)_*=\mathop{\oplus}\limits_{k\geq -d}\fg_k\subset \fv$ with
the given {\bf negative} part.}

{\bf Problem 3}.  {\it Single out {\bf partial prolongs} of
$\fg_-$ in $(\fg_-)_*$. In particular, given not only $\fn$, but
$\fn_0\subset \fder_0\fn$, where the subscript $0$ singles out
derivations that preserve the $\Zee$-grading, we should
automatically have an embedding $\fn_0\subset \fg_0$.

If the inclusion $\fn_0\subset \fg_0$ is a strict one, we wish to
be able to single out $\fn_0$  in $ \fg_0$ as well as to single
out the algebraic prolong $(\fg_-,\fn_0)_*$ --- the maximal
subalgebra of $\fv$ with a given {\bf non-positive} part --- in
$(\fg_-)_*$.

If $\fn_0=\fg_0$ but the component $\fg_1$ forms a reducible
$\fg_0$-module with a submodule $\fh_1$, how to singe out the
maximal subalgebra $(\fg_-\oplus \fg_0\oplus\fh_1)_*\subset \fv$
with a given  \lq\lq beginning part\rq\rq (components of grading
$\le 1$)?

In utmost generality, single out in $\fv$ the maximal subalgebra
$\fh_*=\mathop{\oplus}\limits_{k\geq -d}\fh_k$ with a given {\em
beginning part}
$\fh=\fg_-\oplus\left(\mathop{\oplus}\limits_{0\leq k\leq
K}\fh_k\right)$. Naturally, the beginning part $\fh$ should be
compatible with the bracket, i.e., $[\fh_i,
\fh_j]\subset\fh_{i+j}$ for all $i, j$ such that $i+j\leq K$.}

\vskip 0.2 cm

The components $\fh_k$ with $k>K$ are defined recurrently:
\begin{equation}
\label{2} \fh_k=\{X\in \fg_k\mid [X,\fg_{-1}]\subset
\fh_{k-1}\}. 
\end{equation}
We prove the inclusion $[\fh_k,\fh_l]\subset\fh_{k+l}$ for all
indices by induction on $k+l$ with an appeal to (1); it guarantees
that $\fh_*:=\fg_-\oplus\left(\mathop{\oplus}\limits_{0\leq
k}\fh_k\right)$ is a subalgebra of $\fv$. The Lie algebra
$\fh_*=\mathop{\oplus}\limits_{k\geq -d}\fh_k$ is a generalization
of Cartan prolong.

\begin{Remark} {\it De facto}, for simple Lie algebras over $\Ree$ and $\Cee$,
the number $K$ is always $\leq 1$, but if $\Char \Kee >0$, and for
superalgebras, then $K>1$ is possible. \end{Remark}

{\bf Discussion}. If a Lie group  $N$ with a Lie algebra $\fn$ is
given explicitly, i.e., if we know explicit expressions for the
product of the group elements in some coordinates, then there is
no problem to describe an embedding $\fn\subset \fvect(n)$: the
Lie algebras of left- and right-invariant vector fields on $N$ are
isomorphic to $\fn$. (This is, actually, one of the definitions of
the Lie algebra of $N$.) If the group $N$ is not explicitly given,
then to describe an embedding $\fn\subset \fvect(n)$ is a part of
the problem of recovering the Lie group from its Lie algebra (in
the cases where one can speak about Lie groups). Of course, the
Campbell-Hausdorff formula gives a solution to this problem.
Unfortunately, despite its importance in theoretical discussions,
the Campbell-Hausdorff formula is not convenient in actual
calculations.

For $\Kee=\Ree$ and $\Cee$, another method of constructing an
embedding $\fn\subset \fvect(n)$ and recovering a Lie group from
its Lie algebra is integration of the Maurer-Cartan equations, cf.
\cite{DFN}. Although the algorithm I offer does not use a Lie
group of $\fn$ and is applicable even for the cases where no
analog of a Lie groups can be offered, it is viewing
a given Lie algebra as the Lie algebra of left-invariant vector fields
on a Lie group that gives us a key lead.

For other algorithms for embedding $\fn\subset \fvect(n)$, based
on explicit descriptions of the $\fn$-action in $U(\fn)$, see
\cite{BGLS, M}. Now, let me list reasons that lead to the
algorithm.

{\bf Reason 1.} Let $X_1, \dots, X_n$ be vector fields linearly
independent at each point of an $n$-dimensional (super)domain, and
\begin{equation}\label{3}
[X_i, X_j]=\sum_k c_{ij}^k X_k, \quad c_{ij}^k\in \Kee.
\end{equation}

Let $\omega^1, \dots \omega^n$ be the dual basis of differential
1-forms ($\omega^i(X_j)=\delta^i_j$). Then (a standard exercise)
\begin{equation}\label{4}
d\omega^k= - \frac 12 \sum_{ij} c_{ij}^k \omega^i\wedge \omega^j=
 - \sum_{i<j} c_{ij}^k \omega^i\wedge \omega^j,
\end{equation}
and {\it vise versa}: if the  $1$-forms $\omega^1, \dots \omega^n$
satisfy (\ref{4}) then the dual vector fields $X_1, \dots, X_n$
satisfy (\ref{3}).

Observe that, although in the super setting the expression for
$d\omega$, i.e.,
\begin{equation}\label{dw}
d\omega(X,Y)=X\omega(Y)-Y\omega(X)-\omega([X,Y]),
\end{equation}
acquires some signs, eq. (\ref{4}) is valid for superalgebras as
well: the extra signs in eq. (\ref{dw}) appearing due to super
nature of its constituents do not affect (\ref{4}).

Recall that if the fields $X_i$ form a basis of left-invariant
vector fields on the group $N$, eqs. (\ref{4}) are called the {\it
Maurer-Cartan equations}; in this case, $c_{ij}^k$ are the
structure constants of the Lie algebra $\fn$.

If  $\omega^i=\sum_k V^i_k(x)dx^k$, then eqs. (\ref{4}) can be
expressed as equations for the functions $V^i_k$:
\begin{equation}\label{5}
\partial_j V_i^k-\partial_i V_j^k =\sum_{p,q} c_{pq}^kV_i^pV_j^q.
\end{equation}

{\bf Reason 2.} In the real or complex situation,  eqs. (\ref{5})
are easy to integrate  in \lq\lq nice\rq\rq co\-or\-di\-nates for
any Lie algebra $\fn$ (not only nilpotent). Namely,  introduce
functions
$$
W^i_j(t,x)=tV^i_j(\exp(tx)),\;\text{ where $t\in \Ree$,
$x\in\fn^*$}.
$$
(In other words, we should integrate eqs. (\ref{5}) along
one-parameter subgroups.) As is easy to check, the functions $W$
satisfy ODE
\begin{equation}\label{6}
\frac{dW^i_j}{dt}=\delta^i_j+\sum_{p,q} c_{pq}^i W^p_j x^q
\end{equation}
with the initial condition $W^i_j(0,x)=0$.

Actually, since $\fn$ is  $\Zee$-graded nilpotent, the system
(\ref{5}) is so simple that one can integrate it directly, without
appealing to auxiliary functions $W$, and over any ground field.
This direct solution of (\ref{5}) allows us to construct an
embedding $\fn\tto\fvect(n)$ most suitable for our
purposes\footnote{For example, why even the authors of \cite{BGLS}
were reluctant to use any of the three algorithms presented in
\cite{BGLS}? I tested all of their three algorithms: they work,
although some clarifications (see \cite{Sh2}) are needed. The only
explanation I can deduce from the questions Grozman and Leites
asked me, is the fact that the formulas in \cite{BGLS} are {\it
fixed}, and some of them involve divisions. And what to do if,
say, one wants to avoid division (by 2 or 3) in coefficients?!
Whereas I give the customer a possibility to select the embedding
to taste.}, and find all possible embeddings.

Namely, select a basis $B=\{e_1, \dots, e_n\}$ of $\fn$ compatible
with the grading. This means that its first $n_{1}$ elements form
a basis of $\fn_{-1}$, the next $n_{2}$ elements form a basis of
$\fn_{-2}$, and so on. Let $I_s$ be the set of indices
corresponding to $\fn_{-s}$, and $I=\cup I_s$. Let $c_{ij}^k$ be
the structure constants in this basis:
\begin{equation}\label{basis}
[e_i,e_j]=\sum_{k}c_{ij}^k e_k,
\end{equation}
and $x_1, \dots, x_n$ be the determined by $B$ coordinates of
$\fn^*$, the dual space to $\fn$. The nonstandard $\Zee$-grading
$\fv=\mathop{\oplus}\limits_{k\ge -d}\fv_k$ of $\fvect(n)$
compatible with the $\Zee$-grading of $\fn$ is determined by
setting
\begin{equation}\label{7}
\deg x_i = s\;\text{ for any $i\in I_s$. } 
\end{equation}

Let $X_i\in \fv$ be the image of $e_i$ under our embedding. Then
the value of $X_i$ at 0 is equal to $\partial_{x_i}:=\partial_i$
and the value of the dual form $\omega^i$ at 0 is equal to $dx^i$.
If $i\in I_s$, then  the field $X_i$ and the form $\omega^i$ are
homogeneous of degree $-s$ and $s$ respectively. We have
$$
\renewcommand{\arraystretch}{1.4}
\begin{array}{lll}
\omega^i&=dx^i& \text{for $i\in I_1$};\\
\omega^i&=dx^i+\mathop{\sum}\limits_{j, k\in I_1}a^i_{jk}x^jdx^k& \text{for $i\in I_2$};\\
\omega^i&=dx^i+\mathop{\sum}\limits_{j\in I_1, k\in I_2}a^i_{jk}x^jdx^k+&\\
&\mathop{\sum}\limits_{k\in I_1}\left(\mathop{\sum}\limits_{s,
t\in I_1}a^i_{stk}x^sx^t+
\mathop{\sum}\limits_{s\in I_2}a^i_{sk}x^s\right)dx^k& \text{for $i\in I_3$};\\
\dotfill&\dotfill&\dotfill\\
\end{array}
$$
The grading guarantees automatic fulfillment of a part of
conditions $(\ref{5})$: For example, for $k\in I_1$, all the
functions $V^k_i$ are known: $V^k_i=\delta^k_i$; for $k\in I_2$,
the rhs of $(\ref{5})$ only contains the known functions ($V^k_i$
with $k\in I_1$), and so on.

The system for $a^i_{jk}$ is highly undetermined but if we are
interested in getting some embedding only, we do not need all the
solutions; any solution (the simpler looking, the better) will do.
Then we proceed in the same way with $V^k_i$ for $k\in I_3$, and
so on. The Jacobi identity guarantees the compatibility of the
system.

{\bf Reason 3. How the complete prolongations are singled out.}
Over $\Ree$, the connected simply connected Lie group $N$ with Lie
algebra $\fn$,  left-invariant forms $\omega^i$, where $i\in I$,
and the structure constants  $c_{ij}^k$ given by  (\ref{4})
possess a universal property (\cite{St}):

{\sl Let $M$ be a smooth manifold with a collection of linearly
independent at each point differential 1-forms $\alpha^i$,
satisfying (\ref{4}) with the same constants $c_{ij}^k$. Then, for
every point $x\in M$, there exists its neighborhood $U$ and a
diffeomorphism $f:U\tto N$ such that
$$
\alpha^i=f^*(\omega^i).
$$
Any two such diffeomorphisms differ by a translation.}

Hence, as soon as we have found forms $\omega^i$ satisfying
(\ref{4}), we can think of them as of left-invariant forms of the
group $N$ and of the dual vector fields $X_i$ as of left-invariant
vector fields, $\fg_-= Span\{X_1, \dots, X_n\}\subset \fvect(N)$.

Let $Y_1, \dots, Y_n$ be the right-invariant vector fields, such
that
$$
X_i(e)=Y_i(e),
$$
and $\theta^1, \dots, \theta^n$ be the dual right-invariant
1-forms.

Clearly, both  $\{X_i\}_{i\in I}$ and $\{Y_i\}_{i\in I}$ span Lie
subalgebras of $\fv_-$.

Let us define a right-invariant distribution $\cD$ on $N$ such
that $D(e)=\fn_{-1}$. Clearly, $\cD$ is singled out by the system
of equations for  $X\in\fvect(n)$:
\begin{equation}\label{8}
\theta^{i}(X)=0 \text{ for any } i\in I_2\cup I_3\cup \dots \cup
I_d.
\end{equation}
Since left- and right-invariant vector fields on a Lie group
always commute with each other, each $X_j$ preserves $\cD$ and
hence the Lie algebra $\fg_-$ preserves $ \cD$. Moreover, since
$\fn$ is $\Zee$-graded, it follows that the fact \lq\lq$X\in\fv_-$
preserves $\cD$\rq\rq  is equivalent to the fact \lq\lq$X$
commutes with all $Y_i$, where $i\in I_1$\rq\rq, and hence with
all $Y_i$, where $i\in I$, since $\fn_{-1}$ generates $\fn$.

Thus, $\fg_-$ is characterized as the maximal subalgebra of
$\fv_-$ preserving $\cD$. But then the complete prolongation of
$\fg_-$ is the maximal subalgebra of $\fvect(n)$ preserving $\cD$.

Of course we can reformulate all this without appealing to $N$.
All we need is $\mathfrak{cent}_{\fv_-}(\fg_-)$, the centralizer
of $\fg_-$ in $\fv_-$. It is also clear that, having represented
$Y\in\fv_{-s}$ as a sum of homogeneous components in the standard
grading ($\deg x^i=1$):
$$
Y=\sum_{p=-1}^{d-s-1} Y_{(p)},
$$
we see that for the fields that vanish at the origin (for them,
$Y_{(-1)}=0$), the lowest component of $[X_i,Y]$ coincides with
the bracket of the lowest component of $Y$ with $\partial_i$, and
therefore is nonzero.

The other way round, for any $Y$ such that $Y_{(-1)}\ne 0$ the
equations $[X_i,Y]=0$, where $i=1, \dots, n$, enable us to
uniquely recover, consecutively, all the components $Y_{(p)}$ for
$p\ge 0$ starting with $Y_{(-1)}$ using the recurrence:
\begin{equation}\label{9}
[\partial_i, Y_{(p)}]= -
\sum_{s=-1}^{p-1}[\left(X_i\right)_{(p-1-s)},Y_{(s)}]\quad
\text{for }\; i=1, \dots, n.
\end{equation}
Let $Y_i\in\mathfrak{cent}_{\fv_-}(\fg_-)$ be such that
$\left(Y_i\right)_{(-1)}=\partial_i$. Then
\begin{equation}
\label{er}
\renewcommand{\arraystretch}{1.4}
\begin{array}{l}
[Y_i,Y_j]_{(-1)}=[\partial_i,
\left(Y_j\right)_{(0)}]+[\left(Y_i\right)_{(0)},
\partial_j]= -[\left(X_i\right)_{(0)},
\partial_j]-[\partial_i, \left(X_j\right)_{(0)}]=\\
-[X_i,X_j]_{(-1)}=-\sum_k c_{ij}^k \left(X_k\right)_{(-1)}=-\sum_k
c_{ij}^k \left(Y_k\right)_{(-1)},
\end{array}
\end{equation}
and, since the fields from $\mathfrak{cent}_{\fv_-}(\fg_-)$ are
uniquely determined by their $(-1)$st components, we get:
$$
[Y_i,Y_j]= -\sum_k c_{ij}^k Y_k,
$$
i.e., $\mathfrak{cent}_{\fv_-}(\fg_-)$ is isomorphic to $\fn$.

Let the $\theta^i$ constitute a basis of  1-forms dual to the
$\{Y_i\}_{i\in I}$ (i.e., $\theta^i(Y_j)=\delta^i_j$). Then any
vector field $X\in \fvect(n)$ is of the form
$$
X=\sum_i \theta^i(X)Y_i.
$$
Since $[X_i, Y_j]=0$ for any $i,j=1, \dots, n$, we have
\begin{equation}
\label{10}
\theta^i([X_j,X])=X_j(\theta^i(X)) .
\end{equation}

Now let us consider  the distribution $\cD$ defined by  (\ref{8}).
As we have already observed, $\fg_-$ is characterized as the
maximal subalgebra of $\fv_-$ preserving $\cD$.

Observe, first of all, that although \lq\lq any $X_i$ preserves
any form $\theta^j$\rq\rq, the condition in quotation marks does
not survive the operation (\ref{2}) of complete prolongation
whereas the condition \lq\lq preserve $\cD$\rq\rq\ is not so
strong and survives it.

Indeed, a field $X\in\fv$ preserves $\cD$ if and only if
\begin{equation}
\label{11} \theta^k([X,Y_i])=0 \text{ for any $i=1,\dots, n_1$,
and any $k>n_1$}.
\end{equation}
Let  (\ref{11}) be valid for any $X\in \fg_{s-1}$. Then,  due to
(\ref{2}), $X\in \fg_{s}$ if and only if
\begin{equation}
\label{12}
\theta^k([[X_j,X],Y_i])=X_j\theta^k([X,Y_i])=0\; \text{ for any
$i,j=1,\dots, n_1$, and $ k>n_1$. }
\end{equation}
(We have taken (\ref{10}) into account.)

Finally, since  $\fn_{-1}$ generates the algebra $\fn$, (\ref{12})
is equivalent to
\begin{equation}
\label{12'}
\partial_j(\theta^k([X,Y_i])=0 \;\text{ for all $j=1, \dots, n$}.
\end{equation}

But if $k\in I_l$ ($l\ge 2$), then  $\theta^k([X,Y_i])$ is a
homogeneous (in our nonstandard grading) polynomial of degree
$s-1+l\ge s+1\ge 1$, and hence (\ref{12'}) is equivalent to
\begin{equation}
\label{12''} \theta^k([X,Y_i])=0\; \text{ for any $i=1,\dots,
n_1$, and $ k>n_1$, }
\end{equation}
and hence $X$ preserves $\cD$.

Let us rewrite the system (\ref{11}) for coordinates of $X$ more
explicitly:
\begin{equation}
\label{13}
\renewcommand{\arraystretch}{1.4}
\begin{array}{l}
Y_i(\theta^k(X))-\sum_j(-1)^{p(Y_i)p(\theta^j(X))}c_{ij}^k\theta^j(X)=0\\
\text{for any $i=1, \dots, n_1$, and $k=n_1+1, \dots,
n$}.
\end{array}
\end{equation}

Since $\fg_-$ is $\Zee$-graded, eqs. (\ref{13}) are of a particular
form. Let
$$
F_i=\theta^{n-n_d+i}(X), \; \text{where }\; i=1,\dots, n_d,
$$
be the coordinates of a vector field $X$ lying in the component
$\fg_{-d}$ of maximal depth.

If the functions $F_i$ are given, then eqs. (\ref{13}), where $k\in
I_{d-1}$, constitute a system of linear (not differential)
equations for the coordinates $\theta^j(X)$ corresponding to the
component $\fg_{-d+1}$, and if this component does not contain
central elements of the whole algebra $\fg_-$, then all the
coordinates of the level $-d+1$ enter the system. After all these
coordinates are determined, eqs. (\ref{13}), where now $k\in I_{d-2}$, become a
system of linear equations for coordinates on the next level, $-d+2$,
and so on.

Therefore, the $F_i$ are generating functions for $X$. In the
general case, one should take for generating functions the
functions corresponding to all central  basis  elements of
$\fg_-$.

\medskip

Now, we are able to formulate the algorithm  for the first two of
our problems.

\section{The algorithm: Solving Problems 1 and 2}
$\bullet$ In $\fn$, take a basis $B$ compatible with the grading
and compute the corresponding structure constants $c_{ij}^k$.

$\bullet$ Seek the basis of 1-forms $\{\omega^i\}_{i\in I}$
satisfying (\ref{4}), i.e., solve system (\ref{5})  upwards, i.e.,
starting with degree $1$ and proceeding up to degree $d$.

$\bullet$ Seek the dual basis of vector fields $\{X_i\}_{i\in I}$
upwards, i.e., starting with degree $-d$ and proceeding up to
degree $-1$. The fields $\{X_i\}_{i\in I}$ determine an embedding
of $\fn$ into $\fvect(n)$.

\underline{Problem 1 is solved}.

$\bullet$ Seek a basis $\{Y_i\}_{i\in I}$ of
$\mathfrak{cent}_{\fv_-}(\fg_-)$  in $\fv_-$ by means of (\ref{9})
and the  dual basis of 1-forms $\{\theta^i\}_{i\in I}$.

$\bullet$ To find the component $\fg_s$ of the complete prolongation
of $\fg_-$, we seek the field $X\in \fg_s$ in the form $X=\sum
\theta^i(X)Y_i$. For this, we express each of the $n_d$ generating
functions $F^i=\theta^{n-n_d+i}(X)$ as a sum of monomials of
degree $d+s$ (in the nonstandard grading) with undetermined
coefficients and solve the system (\ref{13})
of linear homogeneous equations for these coefficients.

For debugging, we compare, for $s<0$, the fields thus obtained
with the $X_i$.

\underline{Problem 2 is solved}.

{\bf Example.} Consider the exceptional Lie algebra
$\fg(2)$\footnote{We denote the exceptional Lie algebras in the
same way as the serial ones, like $\fsl(n)$; we thus avoid
confusing $\fg(2)$ with the second component $\fg_2$ of a
$\Zee$-graded Lie algebra $\fg$.} in its $\Zee$-grading of depth
3, as in \cite{C, Y}. In what follows,
$$
x^{(k)}\;\text{ denotes $\begin{cases}\frac{x^k}{k!}&\text{over
$\Ree$ or $\Cee$}\\
u^{(k)}\text{(the divided power)}&\text{in characteristic $p>0$}.
\end{cases}$}
$$
Then (recall that $\fn$ is the given abstract algebra whose image
in the Lie algebra of vector fields is designated by $\fg$)
$$ \fn=\fn_{-3}\oplus
\fn_{-2}\oplus \fn_{-1},\;\text{ where $\dim\fn_{-1}=2$,
$\dim\fn_{-2}=1$, $\dim\fn_{-3}=2$}.
$$
Let us see how the algorithm works for the embedding
$f(\fn)=\fg_-\subset\fvect(5)$.

{\bf 1.} A basis compatible with the $\Zee$-grading and structure
constants are of the form:
$$
\renewcommand{\arraystretch}{1.4}
\begin{array}{l}
\fn_{-1}=\Span(e_1, e_2),\;\; \fn_{-2}=\Span(e_3),\;\;
\fn_{-3}=\Span(e_4,
e_5);\\
{} [e_1,e_2]=e_3, \quad [e_1, e_3]= e_4, \quad [e_2,e_3]=e_5;\\
c_{ij}^k=0 \text{ for } k=1,2; \quad c_{12}^3=-c_{21}^3=1,
c_{ij}^3=0 \text{ otherwise };\\
 c_{13}^4=-c_{31}^4=1, \qquad c_{ij}^4=0
\text{ otherwise;}\\
 c_{23}^5=-c_{32}^5=1, \qquad c_{ij}^5=0 \text{ otherwise. }
\end{array}
$$

{\bf 2.} We have
$$
\renewcommand{\arraystretch}{1.4}
\begin{array}{l}
\omega^1=dx^1, \quad \omega^2=dx^2\Longrightarrow
V^k_i=\delta^k_i\;\text{ for }k=1, 2\\
\omega^3=dx^3+\mathop{\sum}\limits_{i,
j=1}^2a^3_{ij}x^idx^j\Longrightarrow
V^3_4=V^3_5=0;\;\; V^3_3=1;\\
V^3_1=a^3_{11}x^1+a^3_{21}x^2;\quad V^3_2=a^3_{12}x^1+a^3_{22}x^2.\\
\end{array}
$$
Eqs. $(\ref{5})$ give one non-trivial relation on the $V^3_i$:
$$
\partial_2V^3_1-\partial_1V^3_2=V^1_1V^2_2-V^2_1V^1_2=1,
$$
or, equivalently,
\begin{equation}
\label{3*}
a^3_{21}-a^3_{12}=1.
\end{equation}
Select a solution which seems to be a simplest one:
$$
a^3_{11}=a^3_{22}=a^3_{12}=0, \quad a^3_{21}=1.
$$
(In canonical coordinates of first kind, $a^3_{11}=a^3_{22}=0,\;
a^3_{12}=a^3_{21}=\frac 12$. These are most symmetric coordinates.
We wish, however, to evade division if possible.) Thus,
$V^3_1=x^2, \; V^3_2=0$, and hence
$$
\omega^3=dx^3+x^2dx^1.
$$
Further, for $k=4,5$,
$$
\omega^k=dx^k+\sum_{j=1}^3V^k_jdx^j,
$$
where
$$
\renewcommand{\arraystretch}{1.4}
\begin{array}{l}
V^k_3=a^k_1x^1+a^k_2x^2, \\
V^k_j=\alpha^k_j (x^1)^{(2)} + \beta^k_j x^1x^2 +\gamma^k_j
(x^2)^{(2)}+\eps^k_j x^3,\quad j=1,2. \end{array}
$$
Eqs. $(\ref{5})$ give three nontrivial relation for each function $V^k_j$,
where $k=4,5$, $j=1,2,3$:
$$
\renewcommand{\arraystretch}{1.4}
\begin{array}{l}
\partial_2V^4_1-\partial_1V^4_2=0; \quad
\partial_3V^4_1-\partial_1V^4_3=1; \quad
\partial_3V^4_2-\partial_2V^4_3=0;\\
\partial_2V^5_1-\partial_1V^5_2=-x^2; \quad
\partial_3V^5_1-\partial_1V^5_3=0; \quad
\partial_3V^5_2-\partial_2V^5_3=1, \end{array}
$$
or, in terms of coefficients:
$$
\renewcommand{\arraystretch}{1.4}
\begin{array}{l}
\beta^4_1x^1+\gamma^4_1x^2=\alpha^4_2x^1+\beta^4_2x^2; \quad
\eps^4_1-a^4_1=1; \quad \eps^4_2-a^4_2=0;\\
\beta^5_1x^1-\gamma^5_1x^2-\alpha^5_2x^1-\beta^5_2x^2=-x^2; \quad
\eps^5_1-a^5_1=1; \quad \eps^5_2-a^5_2=1. \end{array}
$$
Select a simpler looking solution:
$$
\omega^4=dx^4-x^1dx^3,\quad \omega^5=dx^5-x^2dx^3-(x^2)^{(2)}dx^1.
$$
Finally,
$$
\renewcommand{\arraystretch}{1.4}
\begin{array}{l}
\omega^1=dx^1, \quad \omega^2=dx^2,\\
\omega^3=dx^3+x^2dx^1,\\
\omega^4=dx^4-x^1dx^3,\\
\omega^5=dx^5-x^2dx^3-(x^2)^{(2)}dx^1.\end{array}
$$

{\bf 3.} Now seek the dual fields $X_i$:
$$
\renewcommand{\arraystretch}{1.4}
\begin{array}{l}
X_5=\partial_5, \quad\quad X_4=\partial_4,\\
X_3=\partial_3+x^1\partial_4+x^2\partial_5,\\
X_2=\partial_2, \quad\quad
X_1=\partial_1-x^2\partial_3-x^1x^2\partial_4-(x^2)^{(2)}\partial_5.
\end{array}
$$
We get $\fg_-=\Span\{X_1, \dots, X_5\}$.

{\bf 4.} Now we seek homogeneous fields $Y_i=\partial_i+\dots$,
commuting with all the $X_j$. Since the brackets with $X_2, X_4,
X_5$ vanish, the coordinates of the $Y_i$ can only depend on $x^1$
and $x^3$. Therefore
$$
\renewcommand{\arraystretch}{1.4}
\begin{array}{l}
Y_4=\partial_4, \quad Y_5=\partial_5;\\
Y_3=\partial_3+ ax^1\partial_4+bx^1\partial_5,\end{array}
$$
and $[X_1,Y_3]=0$ implies that $Y_3=\partial_3$.

Finally, for $i=1,2$, we have
$$
Y_i=\partial_i+\alpha_ix^1\partial_3+\sum_{j=4}^5
(\beta_i^j(x^1)^{(2)}+\gamma_i^jx^3)\partial_j.
$$
Bracketing $Y_i$ with $X_1$ and $X_3$, we get
\begin{equation}
\label{kill}
\renewcommand{\arraystretch}{1.4}
\begin{array}{l}
Y_1=\partial_1+x^3\partial_4,\\
Y_2=\partial_2-x^1\partial_3-(x^1)^{(2)}\partial_4+x^3\partial_5.
\end{array}
\end{equation}
It only remains to find the forms $\theta^i$ left-dual to $Y_i$.
The routine computations yield: $\theta^i=dx^i$ for $i=1,2$ and
\begin{equation}
\label{theta}
\renewcommand{\arraystretch}{1.4}
\begin{array}{l}
\theta^3=dx^3+x^1dx^2;\\
\theta^4=dx^4-x^3dx^1+(x^1)^{(2)}dx^2; \\
\theta^5=dx^5-x^3dx^2.\end{array}
\end{equation}

{\bf 5.} Now we seek all the vector fields $X$ preserving
$\cD=\Span\{Y_1, Y_2\}$, or, which is the same, all the fields
that belong to the complete prolongation of $\fg_-$. Let $X=\sum
f^iY_i$, where $f^i=\theta^i(X)$. To find the $f^i$, we solve eqs.
(\ref{13}). In our case they are:
$$
Y_1(f^4)=f^3,\; Y_1(f^5)=0,\; Y_2(f^4)=0,\; Y_2(f^5)=f^3,\;
Y_1(f^3)=f^2,\; Y_2(f^3)=-f^1.
$$
We see that $X$ is completely determined by the functions $f^4$
and $f^5$ which must satisfy the three relations:
\begin{equation}
\label{17}
Y_1(f^5)=0,\; Y_2(f^4)=0,\; Y_1(f^4)=Y_2(f^5).
\end{equation}

For  control, let us look what are the corresponding fields in the
component $\fv_{-2}$. In this case, both $f^4$ and $f^5$ should be
of degree 1 in our grading, i.e., must be of the form
$f^i=a^ix^1+b^i x^2$ for $i=4, 5$. Then $Y_1(f^i)=a^i$,
$Y_2(f^i)=b^i$, and hence eqs. (\ref{17}) mean that
$$
f^4=ax^1, \quad f^5=ax^2 \Longrightarrow f^3=a, f^1=f^2=0.
$$
Therefore, any field preserving $\cal D$ and lying in $\fv_{-2}$
is proportional to
$$
X=Y_3+x^1Y_4+x^2Y_5=\partial_3+ x^1\partial_4+ x^2\partial_5=X_3,
$$
as should be.

We similarly check that, in $\fv_{-1}$, our equations single out
precisely the subspace spanned by $X_1$ and $X_2$.

Now, let us compute $\fg_0$. Its generating functions must be of
degree 3 in our grading, i.e., of the form (for $i=4,5$)
$$
f^i=a_1^i(x^1)^{(3)}+ a_2^i(x^1)^{(2)}x^2+ a_3^ix^1(x^2)^{(2)} +
a_4^i(x^2)^{(3)} + b_1^ix^1x^3+ b_2^ix^2x^3+ c_1^ix^4 + c_2^ix^5.
$$
Then
$$
\renewcommand{\arraystretch}{1.4}
\begin{array}{l}
Y_1(f^i)= a_1^i(x^1)^{(2)}+ a_2^ix^1x^2+ a_3^i(x^2)^{(2)} +
(b_1^i+c_1^i)x^3;\\
 Y_2(f^i)=
(a_2^i-2b_1^i-c_1^i)(x^1)^{(2)}+ (a_3^i-b_2^i)x^1x^2+
a_4^i(x^2)^{(2)} +(b_2^i+c_2^i)x^3.\end{array}
$$
In this case, eqs. (\ref{17})  take the following form:
$$
\left\{\renewcommand{\arraystretch}{1.4}
\begin{array}{l}
a_1^5=a_2^5=a_3^5=0,\\
 b_1^5+ c_1^5=0, \\
a_4^4=0,\\
 b_2^4+ c_2^4=0, \\
  a_3^4- b_2^4=0,\\
 a_2^4-2 b_1^4+c_1^4=0,\\
a_1^4=a_2^5- 2 b_1^5-c_1^5,\\
 a_2^4=a_3^5- b_2^5,\\
a_3^4=a_4^5,\\
 b_1^4+c_1^4= b_2^5+ c_2^5.\end{array}\right.
$$
The solution to this system is:
$$
\renewcommand{\arraystretch}{1.4}
\begin{array}{l}
a_1^4=-b_1^5=c_1^5=\alpha,\\ a_2^4=-b_2^5=\beta,\\
a_3^4=a_4^5=b_2^4=-c_2^4=\gamma,\\ a_4^4=a_1^5=a_2^5=a_3^5=0,\\
b_1^4=\delta,\\ c_1^4=\beta-2\delta,\\
 c_2^5=2\beta-\delta.\end{array}
$$
Hence
$$
\renewcommand{\arraystretch}{1.4}
\begin{array}{l}
f^4= \alpha(x^1)^{(3)}+ \beta(x^1)^{(2)}x^2+ \gamma x^1(x^2)^{(2)}
+
 \delta x^1x^3+ \gamma x^2x^3+ (\beta-2\delta)x^4 -\gamma
 x^5,\\
f^5= \gamma (x^2)^3 - \alpha
 x^1x^3-\beta x^2x^3+ \alpha x^4 +(2\beta-\delta)
 x^5,\\
f^3= \alpha(x^1)^{(2)}+ \beta x^1x^2+ \gamma (x^2)^{(2)} +
(\beta-\delta) x^3,\\
f^2=Y_1(f^3)=\alpha
x^1+\beta x^2,\\
 f^1=-Y_2(f^3)=-\delta x^1-
\gamma x^2.\end{array}
$$
For a basis of $\fg_0$ we take the vectors $X_\alpha, X_\beta,
X_\gamma, X_\delta$ corresponding to the only one non-zero
parameter (for example, $X_\alpha$ corresponds to $\alpha=1$,
$\beta=\gamma=\delta=0$ and so on):
$$
\renewcommand{\arraystretch}{1.4}
\begin{array}{ll}
X_\alpha &=x^1Y_2+  (x^1)^{(2)}Y_3 +(x^1)^{(3)}Y_4+
(-x^1x^3+x^4)Y_5=\\
&x^1\partial_2+x^4\partial_5-(x^1)^{(2)}\partial_3-2(x^1)^{(2)}\partial_4,\\
X_\beta&=x^2\partial_2+x^3\partial_3+x^4\partial_4+2x^5\partial_5,\\
X_\gamma& =- x^2\partial_1-x^5\partial_4
+(x^2)^{(2)}\partial_3+x^1(x^2)^{(2)}\partial_4+
(x^2)^{(3)}\partial_5, \\
X_\delta&= -x^1\partial_1-x^3\partial_3 - 2 x^4\partial_4
-x^5\partial_5.\end{array}
$$
For $\alpha=\gamma=0$, $\delta=-\beta=1$, we get the grading
operator
$$
X= -x^1\partial_1-x^2\partial_2 -2x^3\partial_3- 3x^4\partial_4 -3
x^5\partial_5.
$$
The higher components can be calculated in a similar way.

{\bf Interpretation}. There are three realizations of $\fg=\fg(2)$
as a Lie algebra that preserves a non-integrable distribution on
$\fg_-$ related with the three (incompressible) $\Zee$-gradings of
$\fg$: with one or both coroots of degree 1. Above we considered
the grading $(1, 0)$; Cartan used it to give the first
interpretation of $\fg(2)$, then recently discovered by Killing,
see \cite{C}.\footnote{Cartan also considered the grading $(1,
1)$, see \cite{Y}; Cartan used it to study Hilbert's equation
$f'=(g'')^2$. The grading $(0, 1)$ is considered in \cite{La2} but
we could not follow his calculations and intend to redo them.}

In this realization (by fields $X_i$) $\fg=\fg(2)$ preserves the
distribution in the tangent bundle on $\fg_-$ given by the system
of Pfaff equations for vector fields $X$
$$
\renewcommand{\arraystretch}{1.4}
\begin{array}{l}
\theta^3(X)=0;\qquad \theta^4(X)=0; \qquad
\theta^5(X)=0.\end{array}
$$
Equivalently, but a bit more economically, we can describe
$\fg=\fg(2)$ as preserving the codistribution in the cotangent
bundle on $\fg_-$ given by the vectors $(\ref{kill})$, i.e., as
the following system of equations for 1-forms $\alpha$:
$$
\alpha(Y_1)=0;\qquad \alpha(Y_2)=0.
$$

Obviously, description in terms of codistributions is sometimes
shorter: any distribution of codimension $r$ requires for its
description $r$ Pfaff equations, whereas the dual codistribution
requires $n-r$ equations.

One can similarly describe the remaining realizations of $\fg(2)$
corresponding to the other $\Zee$-gradings, various realizations
of $\ff(4)$ and $\fe(6)-\fe(8)$ and of exceptional Lie
superalgebras, as well as Lie algebras over fields of
characteristic $p$. There seemed to be no need to consider
nonintegrable distributions associated with various
$\Zee$-gradings of {\it non-exceptional} Lie algebras (their usual
description as preserving volume or a nondegenerate form seems to
be sufficiently clear); Cartan himself, though understood
importance of description of Lie algebras in terms of
distributions, only considered one or two $\Zee$-gradings and
related distributions of exceptional Lie algebras and none for
non-exceptional. If, however, we apply the algorithm presented
here to $\fg(2)$, $\fo(7)$, $\fsp(4)$ and $\fsp(10)$ in
characteristic $p=2, 3$ or 5, we elucidate the meaning of some of
the simple Lie algebras specific to $p=2, 3, 5$ and, with luck and
in the absence of classification, distinguish new examples, as in
\cite{GL}.

Other gradings of other algebras are now being under
consideration.

\section{How to single out partial prolongs: Solving Problem 3}
Thus, we have described the {\bf complete} prolong of the Lie
(super)algebra $\fg_-$, i.e., as we have already observed, the
maximal subalgebra $\fg=(\fg_-)_*\subset \fv$ with a given
negative part. Let us consider now a subspace
$\fh=\fg_-\oplus\left(\mathop{\oplus}\limits_{0\leq k\leq
K}\fh_k\right)\subset \fg$, closed with respect to the bracket
within limits of its degrees, i.e., such that
$[\fh_i,\fh_j]\subset \fh_{i+j}$ whenever $i+j\le K$. Let us
describe the {\bf partial} prolong
$\fh_*=\mathop{\oplus}\limits_{k\geq -d}\fh_k\subset \fg$ of the
subspace $\fh$, i.e., the maximal subalgebra of $\fg$ with the
given {\it beginning part} $\fh$. The components $\fh_k$ with
$k>K$ are singled out by condition (\ref{2}). Here by description
we mean a way to single out $\fh$ in $\fg$ by a system of
differential equations.

\begin{Remark} Observe that the Cartan prolong $(\fg_-,
\fg_0)_*$ (where $\fg_-$ is commutative, the depth is $d=1$, and
$\fg_0\subset \fgl(n)$) is a particular case of the above
construction with $\fg=\fvect(n)$, and $\fh=\fg_-\oplus \fg_0$.
For examples of descriptions of Cartan prolongations by means of
differential equations, see\cite{Sh14}, \cite{ShFA}, and
\cite{ShP}. \end{Remark}

The homogeneous component $\fh_m$ of $\fh$ is said to be {\it
defining}, if $\fh_k=\fg_k$ for all $k<m$ but $\fh_m\ne \fg_m$.
Let us consider an algorithm of description of $\fh_*$ in the case
where the defining component is of the maximal degree --- $\fh_K$.
The case of defining component of smaller degree $m<K$ can be
reduced to our case; indeed, we first describe the partial prolong
$\fh^m_*=\left(\mathop{\oplus}\limits_{0\leq k\leq
m}\fh_k\right)_*$, compare the components  $\fh_k$, where $m<k\le
K$ with the corresponding components of this prolong $\fh^m_*$,
find out the new defining component, if any, and so on.

Thus, let $Z_1, \dots, Z_{\dim\fh_K}$ be a basis of the defining
component $\fh_K\subset \fg_K$.

The first thing to do is to single out the subspace $\fh_K$ in
$\fg_K$ by means of a system of linear (algebraic) equations
(i.e., find out a  basis of the annihilator of $\fh_K$ in
$(\fg_K)^*$, or, equivalently, find out the fundamental system of
solutions $\alpha^1, \dots, \alpha^r$ of the system of equations
for an unknown 1-form $\alpha\in (\fg_K)^*$:
\begin{equation}
\label{alpha} \alpha(Z_i)=0 \text{ for all } i=1, \dots,
\dim\fh_K.
\end{equation}
The subspace $\fh_K$ is then singled out by a system of
homogeneous linear equations for an unknown vector field
$X\in\fg_K$:
\begin{equation}
\label{eqhK} \alpha^i(X)=0 \text{ for all } i=1, \dots, r.
\end{equation}
Observe now that in $\fg_K$ there is a convenient for us basis
consisting of the fields of the form $fY_j$, where $f$ is a
monomial of degree $K+s$ if $j\in I_s$. Accordingly, the dual
basis consists of the elements of degree $-K$ and of the form
$$
A_{i_1, \dots, i_t}^j= S(Y_{i_1}\dots Y_{i_t})\theta^j,
$$
and the forms $\alpha^l$ can be expressed in this basis as
\begin{equation}
\label{alphal} \alpha^l=\sum a_j^{l;i_1, \dots, i_t}A_{i_1, \dots,
i_t}^j.
\end{equation}
Substituting (\ref{alphal}) into (\ref{eqhK}) we get a system of
homogeneous linear differential equations with constant
coefficients for the coordinates of the vector field
$X=\theta^i(X)Y_i\in \fh_K$:
\begin{equation}
\label{difeq} \sum a_j^{l;i_1, \dots, i_t}S(Y_{i_1}\dots
Y_{i_t})\theta^j(X) \text{ for all } l=1, \dots, r.
\end{equation}

Observe now that, for Lie algebras, equation (\ref{difeq}) survive
prolongation procedure (\ref{2}). Indeed, for $k>K$, by the
induction hypothesis $X\in \fh_k$ if and only if the brackets
$[X_i,X]$ satisfy (\ref{difeq}) for any $i=1,\dots, n_1$. Set
$$
f^l= \sum a_j^{l;i_1, \dots, i_t}S(Y_{i_1}\dots
Y_{i_t})\theta^j(X).
$$
Since all the $X_i$ commute with all the $Y_j$, the system
(\ref{difeq}) for the brackets $[X_i,X]$ is equivalent to the
system
\begin{equation}
\label{difeq2} X_i(f^l)=0 \text{ for all $i=1, \dots, n_1$ and
$l=1, \dots, r,$}
\end{equation}
which thanks to (\ref{1}) is, in its turn, equivalent to the
system
$$
\partial_i(f^l)=0 \text{ for all $i=1, \dots, n$ and $l=1, \dots,
r.$}
$$
This implies that $f^l=\text{const}$ for all $l=1,\dots, r$. Since
the functions  $f^l$ are homogeneous polynomials of degree
$k-K>0$, it follows that $f^l=0$. Hence, $X\in \fh_k$ if and only
if $X$ satisfies system (\ref{difeq}).

In super case the fields $X_i$ and $Y_j$ supercommute, not
commute, and this does not allow us, generally speaking, break out
the $X_i$ and pass from system (\ref{difeq}) for the brackets to
the system (\ref{difeq2}). There is, however, a simple and
well-know consideration that saves us. Recall that $p$ denotes the
parity function and $\Pty$ is the parity operator, i.e.,
$$
\Pty(x)=(-1)^{p(x)}x.
$$

\begin{Lemma}
Let $X, Y\in \End V$ supercommute and $p(Y)=\od$. Then $X$ and
$\hat Y=Y\Pty$ commute (in the usual sense), i.e., $X\hat Y=\hat Y
X$.\end{Lemma}

Indeed,
$$
\renewcommand{\arraystretch}{1.4}
\begin{array}{ll}
X\hat
Y(v)&=XY\Pty(v)=(-1)^{p(v)}XY(v)=\\
&(-1)^{p(v)}(-1)^{p(X)p(Y)}YX(v)= (-1)^{p(v)+p(X)}YX(v)= \hat
YX(v).\end{array}
$$
Therefore, in the super case, the system (\ref{difeq}) should be
written with operators $\hat Y_i$ instead of $Y_i$ (if
$p(Y_i)=\ev$, we set $\hat Y_i= Y_i$).

Finally, if $d>1$, then any field $X\in \fg$ is completely
determined by its generating functions $F^i$. Therefore, it
suffices to write equations (\ref{difeq}) for the generating
functions only.

{\bf Examples: \underline{Depth 1}}.  Let $\fg_-=\fg_{-1}$ be
commutative, hence $\fg=(\fg_-)_*=\fvect(n)$ in the standard
$\Zee$-grading (the degree of each indeterminate is equal to 1).
Let $\fg_0=\fgl(n)=\fvect(n)_0$. The degree 1 component
$\fvect(n)_1$ consists, as is well-known, of the two irreducible
$\fgl(n)$-modules. Let
\begin{equation}
\label{d} X=\sum a^k_{ij}x^ix^j\partial_k:=\sum f^k(x)\partial_k.
\end{equation}
Then these submodules are:
\begin{equation}
\label{w1} \fh_{1(1)}:=\Span\{x^i\sum_j x^j\partial_j\mid i=1,
\dots , n\}
\end{equation}
and
\begin{equation}
\label{w2}
\fh_{1(2)}:=\Span\{X=\sum_{j=1}^n
d^k_{ij}x^ix^j\partial_k \mid d^i_{ii}+\sum d^j_{ij}=0\;\text{ for
}\;i=1, \dots , n.\}
\end{equation}
Let us single out the partial prolongs
$\fh_{*(i)}=(\fg_{-1}\oplus\fgl(n)\oplus\fh_{1(i)})_*$, where
$i=1,2$ in $\fvect(n)$ by means of differential equations on the
functions $f^k(x)$, see (\ref{d}). In this case $X_i=Y_i=\partial_i$.

The conditions on $\fh_{1(2)}$ can be immediately expressed as
$$
\sum_{j}\partial_i\partial_j(f^j)=0\;\text{ for all }\;i=1, \dots
, n
$$
or, equivalently, as
\begin{equation}
\label{ddiv}
\partial_i\left(\sum_{j}\displaystyle\frac{\partial
f^j}{\partial x^j}\right)=0\;\text{ for all }\;i=1, \dots ,
n.
\end{equation}
This is exactly the system (\ref{difeq}) for $\fh_{*(2)}$ which
can be rewritten in a well-known way:
 $$
 \sum_{j}\displaystyle\frac{\partial
f^j}{\partial x^j}=\Div X=\const.
$$
Hence, as is well-known,
$\fh_{*(2)}=\fd\fsvect(n):=\fsvect(n)\subplus\Kee E$, where
$E=\sum x^i\partial_i$.

Now let us consider $\fh_{*(1)}$ (which is, of course, $\fsl(n+1)$
embedded into $\fvect(n)$). Having expressed $X\in \fh_{1(1)}$ as
$$
(c_1(x^1)^2+c_2x^1x^2+\dots+c_nx^1x^n)\partial_1+\dots+
(c_1x^1x^n+c_2x^2x^n+\dots+c_n(x^n)^2)\partial_n
$$
we immediately see that $d^k_{ij}=0$ if $i\neq k$ and $j\neq k$,
and $d^k_{kk}=d^i_{ki}$ for any $i\neq k$. The corresponding
system of differential equations is
\begin{equation}
\label{sl}
\renewcommand{\arraystretch}{1.6}
\begin{array}{l}
\displaystyle\frac{\partial^2f^k}{\partial x^i\partial x^j}=0\;\text{ for $i, j\neq
k$};\\
\displaystyle\frac12\displaystyle\frac{\partial^2f^k}{\left(\partial x^k\right)^2}=
\displaystyle\frac{\partial^2f^i}{\partial x^i\partial x^k}\;\text{ for $i\neq
k$}.
\end{array}
\end{equation}
{\bf Superization}. For superalgebras, as we have seen, one should take compositions
of $Y_i=\partial_i$ with the parity operators, i.e., instead of
the $\partial_i$ we should take operators
$$
\nabla_i(f):=(-1)^{p(f)p(\partial_i)}\partial_i(f).
$$
These
$\nabla_i$ commute (not {\it super}commute) with any operator
$X_j=\partial_j$ from $\fg_{-1}$. The system (\ref{ddiv}) will
take form
$$
\nabla_i\left(\sum_j \nabla_j(f^j)\right)=0 \text{ for all } i=1,
\dots, n,
$$
which yields, nevertheless, the same condition $\Div X =
\text{const}$. (This is one more way to see why the coordinate
expression of divergence in the super case must contain some signs:
eqs. (\ref{ddiv}) do not survive the prolongation procedure
(\ref{2}).)

Having in mind  that
\begin{equation}
\label{sup}
d^k_{ij}=-(-1)^{p(f^k)}\displaystyle\frac{\partial^2f^k}{\partial
x^i\partial x^j}; \;\text{ so $d^k_{kk}=0$ for $x^k$ odd}
\end{equation}
we deduce that the second line in (\ref{sl}) takes the following
form
\begin{equation}
\label{supsl}
\renewcommand{\arraystretch}{1.6}
\begin{array}{l}
\displaystyle\frac12\displaystyle\frac{\partial^2f^k}{\left(\partial x^k\right)^2}=
(-1)^{p(x^i)p(f^i)+1}\displaystyle\frac{\partial^2f^i}{\partial x^i\partial x^k}\;
\text{ for $p(x^i)=\ev$ and $i\neq k$};\\
(-1)^{p(x^j)p(f^j)}\displaystyle\frac{\partial^2f^i}{\partial x^j\partial x^k}=
(-1)^{p(x^i)p(f^i)}\displaystyle\frac{\partial^2f^i}{\partial x^i\partial x^k}\;
\text{ for $p(x^i)=\od$ and $i, j\neq k$}.
\end{array}
\end{equation}

{\bf \underline{Depth $>1$}} We consider several more or less
well-known examples and a new one ($\fk\fas$).

Let
$\fn=\fn_{-1}\oplus\fn_{-2}$ be the Heisenberg Lie algebra: $\dim
\fn_{-1}=2n$, $\dim\fn_{-2}=1$. The complete prolong of $\fn$ is
the Lie algebra  $\fk(2n+1)$ of contact vector fields. Having
embedded $\fn$ into $\fvect(p_1, \dots, p_n, q_1, \dots, q_n; t)$
with the grading $\deg p_i=\deg q_i=1$ for all $i$ and $\deg t=2$
we can take for the $X$-vectors, for example,
$$
X_{q_i}=\partial_{q_i}+p_i\partial_t,\quad X_{p_i}=\partial
_{p_i}-q_i\partial _t;\quad X_t=\partial_t.
$$
Hence $\fg_-=\Span\{X_{p_1}, \dots, X_{p_n}, X_{q_1}, \dots,
X_{q_n}, X_t\}$ and the contact vector fields in consideration
preserve the distribution $\cD$ given by the Pfaff equation
$\alpha(X)=0$ for vector fields $X$, where  $\alpha=dt+\sum_i
(p_idq_i-q_idp_i)$.

The $Y$-vectors in this case are of the form
$$
Y_{q_i}=\partial_{q_i}-p_i\partial_t,\quad Y_{p_i}=\partial
_{p_i}+q_i\partial _t;\quad Y_t=\partial_t.
$$

In this particular example, a contact vector field $K$ is
determined by only one generating function $F$ which is exactly
the coefficient of $Y_t$ in the decomposition of $K$ with respect
to the $Y$-basis and there are no restrictions on the function
$F$. Denoting $F=2f$ and solving eqs. (\ref{13}), we get the
formula for any contact vector field $K_f$:
\begin{equation}
\label{K_f}
\renewcommand{\arraystretch}{1.6}
\begin{array}{ll}
K_f&=2fY_t+\sum_i (-Y_{q_i}(f)Y_{p_i}+Y_{p_i}(f)Y_{q_i})=\\
&(2-E)(f)\partial_t+\frac{\partial f}{\partial
t}E+\sum_i\left(\frac{\partial f}{\partial_{p_i}}\partial
_{q_i}-\frac{\partial f}{\partial_{q_i}}\partial_{p_i}\right),
\end{array}
\end{equation}
where $E=\sum_i(p_i\partial_{p_i}+q_i\partial_{q_i})$. Of course,
this is exactly the standard formula of the contact vector field
with the generating function $f$. Further we use the realization
of $\fk(2n+1)$ in generation functions $f$.

{\underline{If $\fh_0\neq \fk_0$}}, then $\fh_0$ is the defining
component. The component $\fk_0$ is generated by 2nd order
homogeneous polynomials in $p, q, t$. Thus, for a basis $Z$ in
$\fk_0$ we can take monomials $p_ip_j$, $q_iq_j$, $p_iq_j$; $t$
and for a basis $Z^*$ of the dual space we then take operators
$$
Y_{p_i}Y_{p_j}, \;\;Y_{q_i}Y_{q_j}, \;\;Y_{p_i}Y_{q_j} \;\text{for
$i\neq j$};
\;\;\frac12\left(Y_{p_i}Y_{q_i}+Y_{q_i}Y_{p_i}\right);\quad Y_{t}.
$$

To describe the complete prolongation of $\fg_-\oplus \fh_0$, one
should first single out $\fh_0$ in $\fk_0$ in terms of equations
for  functions generating $\fh_0$ in basis $Z$, then rewrite the
equations in terms of $Z^*$.

{\underline{1) If $\fh_0=\fsp(2n)$}}, the generating functions do
not depend on $t$, which means that
$$
Y_t(f)=0.
$$
This equation singles out in $\fk(2n+1)$ the Poisson subalgebra
$\fp\fo(2n)$.

{\underline{2) If $\fh_0\simeq\Cee\ \id$}}, the generating
function is $t$, which means
\begin{equation}\label{cont}
Y_{p_i}Y_{p_j}(f)=0, \;\;Y_{q_i}Y_{q_j}(f)=0,
\;\;Y_{p_i}Y_{q_j}(f)=0 \;\text{for $i\neq j$};
\;\;\left(Y_{p_i}Y_{q_i}+Y_{q_i}Y_{p_i}\right)(f)=0.
\end{equation}

For $i\neq j$, these equations imply
$$
Y_{q_i}Y_{p_i}Y_{p_j}(f)-Y_{p_i}Y_{q_i}Y_{p_j}(f)=Y_tY_{p_j}(f)=Y_{p_j}\left(Y_t(f)\right)=0.
$$
Analogously, $Y_{q_j}\left(Y_t(f)\right)=0$, and hence
$Y_{t}\left(Y_t(f)\right)=0$, i.e., $Y_t(f)=const$ and $f=ct+f_0$
while eqs. (\ref{cont}) imply $\deg_{p,q}f_0\le 1$. Hence the
prolong of $\fg_-\oplus \fh_0$ coincides with $\fg_-\oplus \fh_0$.

{\underline{Let $\fh_0=\fk_0$, and $\fh_1\subset \fk_1$}. As a
$\fk_0$-module, $\fk_1$ decomposes into the direct sum of two
(over $\Cee$; for $\Char \Kee=3$ and in super setting, even over
$\Cee$, the situation is more involved) irreducible submodules,
$W_1$ spanned by cubic monomials in $p$ and $q$, and $W_2$ spanned
by $tp_i$ and $tq_i$. The dual bases of $W_1$ and $W_2$ are given
by order 3 symmetric polynomials in the $Y_{p_i}, Y_{q_i}$, and,
respectively, spanned by $Y_{p_j}Y_t$ and $Y_{q_j}Y_t$.

Hence the subspace $W_1$ is singled out by conditions
$$
Y_{p_j}Y_t(f)=Y_{q_j}Y_t(f)=0 \Longrightarrow Y_t(f)=const.
$$
This equation singles out in $\fk(2n+1)$ the derivation algebra
$\fder(\fp\fo(2n))=\fp\fo(2n)\oplus\Cee K_t$  of the Poisson
algebra.

To single out $W_2$, we have the system
\begin{equation}\label{con2}
\renewcommand{\arraystretch}{1.4}
\begin{array}{l}
Y_{p_i}Y_{p_j}Y_{p_k}(f)=0, \;\;Y_{q_i}Y_{q_j}Y_{q_k}(f)=0,\\
Y_{p_i}Y_{p_j}Y_{q_k}(f)=0, \;\;Y_{p_k}Y_{q_i}Y_{q_j}(f)=0,
\;\;\text{for
$k\neq i, j$};\\
Y_{p_i}\left(Y_{p_j}Y_{q_j}+Y_{q_j}Y_{p_j}\right)(f)=0, \;
Y_{q_i}\left(Y_{p_j}Y_{q_j}+Y_{q_j}Y_{p_j}\right)(f)=0
\;\;\text{for
$i\neq j$}; \\
\left(Y_{p_i}^2Y_{q_i}+Y_{p_i}Y_{q_i}Y_{p_i}+Y_{q_i}Y_{p_i}^2\right)(f)=0,\\
\left(Y_{q_i}^2Y_{p_i}+Y_{q_i}Y_{p_i}Y_{q_i}+Y_{p_i}Y_{q_i}^2\right)(f)=0.
\end{array}
\end{equation}
which implies that $Y_t(f)$ satisfies eqs. (\ref{cont}), and hence
$$
\frac{\partial f}{\partial t}\in \Span\{1;\ p_1, \dots, p_n, \
q_1, \dots, q_n,;\ t\},
$$
whereas
$$
f\in\Span\{ 1;\ p_1, \dots, p_n,\ q_1, \dots, q_n; \ t, \ tp_1,
\dots, tp_n, \ tq_1, \dots, tq_n;\ t^2\}.
$$
Hence the complete prolongation $(\fk_-\oplus \fk_0\oplus W_2)_*$
is isomorphic to $\fsp(2n+2)$.

\underline{ $\fk\fa\fs\subset\fk(1|6)$}. Let $\fn=\fn_{-1}\oplus
\fn_{-2}$ be the Heisenberg Lie superalgebra $\fh\fe\fii(1|6)$:
$\dim \fn_{-1}=0|6$, $\dim\fn_{-2}=1|0$. The complete prolong of
$\fn$ is $\fg=\fk(1|6)$ with $\fg_0=\fc\fo(6)$. The component
$\fg_1$ consists of three irreducible $\fg_0$-modules.

If we consider $\fk(1|6)$ in realization by generating functions
in $t, \theta_1, \dots, \theta_6$, i.e., when
$$
K_f=(2-E)(f)\partial_t+ \frac{\partial f}{\partial
t}E-(-1)^{p(f)}\sum_i \frac{\partial f}{\partial
\theta_i}\partial_{\theta_i}, \quad f\in\Cee[t, \theta_1, \dots,
\theta_6],
$$
where $E=\sum_i \theta_i\partial_{\theta_i}$ and $\{\theta_i,
\theta_j\}_{k.b.}=\delta_{ij}$, then $\fg_1\simeq
t\Lambda(\theta)\oplus \Lambda^3(\theta)=t\Lambda(\theta)\oplus
\fg_1^+\oplus \fg_1^-$ with $\fg_1^{\pm}\subset \Lambda^3(\theta)$
singled out with the help of the Hodge star ${}^*$:
\begin{equation}
\label{g1} \fg_1^{\pm}=\{f\in\Lambda^3(\theta)\mid
f^*=\pm\sqrt{-1}f\}.
\end{equation}
Recall that the Hodge star ${}^*$ is just the Fourier
transformation in odd indeterminates whereas $t$ is considered a
parameter:
$$
\ast: f(\xi, t)\mapsto f^*(\eta, t)=\int \exp(\sum
\eta_i\xi_i)f(\xi, t)\vvol(\xi).
$$

The exceptional simple Lie superalgebra $\fk\fa\fs$ is defined
(\cite{Sh14, CK6}) as a partial  prolong of
$\fh=\mathop{\oplus}\limits_{k=-2}^1 \fh_k$, where
$\fh_k=\fk(1|6)_k$ for $-2\le k\le 0$ and where
$\fh_1=t\Lambda(\theta)\oplus \fg_1^+$. Hence $\fh_1$ is the
defining component.

Then
$$
X_i=\partial_{\theta_i}+\theta_i\partial_{t} \text{ for }\; i=1,
\dots 6,  \quad X_7=\partial_{t},
$$
and
$$
Y_i=\partial_{\theta_i}-\theta_i\partial_{t}\text{ for }\; i=1,
\dots 6,    \quad Y_7=X_7.
$$

Let $I=\{i_1, i_2, i_3\} \subset \{1, \dots, 6\}$ be an ordered
subset of indices, and $I^*=\{j_1, j_2, j_3\}$ the dual subset of
indices (i.e., $\{I,I^*\}$ is an even permutation of $\{1, \dots,
6\}$). Set:
$$
Y_I=Y_{i_1}Y_{i_2}Y_{i_3},  \quad
 Y_{I^*}=Y_{j_1}Y_{j_2}Y_{j_3},
$$
and define $\Delta_{Y_I}: \Cee[t, \theta]\tto \Cee[t, \theta]$ by
the eqs.
$$
\Delta_{Y_I}(f)=(-1)^{p(f)} Y_I(f).
$$
Observe that $\Delta_{Y_I}(t\theta_s)=0$ for any $s=1, \dots, 6$.
Therefore $\fh_1$ can be singled out in $\fk(1|6)_1$ by the
following 10 equations parameterized by partitions $(I,I^*)$ of
$(1, \dots, 6)$ constituting even permutations:
\begin{equation}
\label{K1} (\Delta_{Y_I}-\sqrt{-1}\Delta_{Y_I^*})(f)=0.
\end{equation}
Clearly, (\ref{K1}) is equivalent to
$$
(Y_I-\sqrt{-1}Y_{I^*})(f)=0 .
$$
The solutions of this system span the following subspace of the
space of generating functions:
\begin{equation}
\label{K2}
\renewcommand{\arraystretch}{1.4}
\begin{array}{l}
    f(t)-\sqrt{-1}f'''(t)1^*, \\ 
f_j(t)\theta_j-\sqrt{-1}f''_j(t)\theta_j^*,\\
f_{jk}(t)\theta_j\theta_k-\sqrt{-1}f'_{jk}(t)(\theta_j\theta_k)^*,\\
f_{jkl}(t)\left(\theta_j\theta_k\theta_l-\sqrt{-1}(\theta_j\theta_k\theta_l)^*\right).
\end{array}
\end{equation}
In (\ref{K2}), $j, k, l$ are distinct indices 1 to 6.

The above description is in agreement with equations from
\cite{CK6}.


\begin{thebibliography}{9999}

\bibitem[BGLS]{BGLS}
Burdik, Ch.; Grozman, P.; Leites, D.; Sergeev, A. A construction
of Lie algebras and superalgebras by means of creation and
annihilation operators. I. (Russian. Russian summary) Teoret. Mat.
Fiz. 124 (2000), no. 2, 227--238; translation in Theoret. and
Math. Phys. 124 (2000), no. 2, 1048--1058

\bibitem[C]{C}
Cartan \'E., \"Uber die einfachen Transformationsgrouppen,
Leipziger Berichte (1893), 395--420. Reprinted in: {\em \OE uvres
compl\`{e}tes}. Partie II. (French) [Complete works. Part II]
Alg\`{e}bre, syst\`{e}mes diff\'erentiels et probl\`{e}mes
d'\'equivalence. [Algebra, differential systems and problems of
equivalence] Second edition. \'Editions du Centre National de la
Recherche Scientifique (CNRS), Paris, 1984.

\bibitem[CK6]{CK6}
Cheng, Shun-Jen; Kac, V. A new $N=6$ superconformal algebra. Comm.
Math. Phys. 186 (1997), no. 1, 219--231

\bibitem[DFN]{DFN}
Dubrovin, B. A.; Fomenko, A. T.;
Novikov, S. P. {\em Modern geometry---methods and applications.
Part I. The geometry of surfaces, transformation groups, and
fields}. Second edition. Translated from the Russian by Robert G.
Burns. Graduate Texts in Mathematics, 93. Springer-Verlag, New
York, 1992. xvi+468 pp.

\bibitem[FSh]{FSh}
Fei, Q.-Y., Shen, G.-Y., Universal graded Lie algebras, J. Algebra
152 (1992), 439--453

\bibitem[Gr]{Gr} Grozman P., {\bf SuperLie},
\text{http://www.equaonline.com/math/SuperLie}

\bibitem[GL]{GL}
Grozman P., Leites D., Structures of $G(2)$ type and nonintegrable
distributions in characteristic $p$, arXiv:  math.RT/0509400

\bibitem[La1]{La1}
Larsson T., Structures preserved by consistently graded Lie
superalgebras, arXiv: math-ph/0106004


\bibitem[La2]{La2}
Larsson T., Structures preserved by exceptional Lie algebras,
arXiv: math-ph/0301006

\bibitem[M]{M}
Molotkov V., Explicit realization of induced and coinduced modules
over Lie superalgebras by differential operators, arXiv:
math.RT/0509105

\bibitem[ShFA]{ShFA}
Shchepochkina I., Five simple exceptional Lie superalgebras of
vector fields. (Russian) Funktsional. Anal. i Prilozhen. 33
(1999), no. 3, 59--72, 96; translation in Funct. Anal. Appl. 33
(1999), no. 3, 208--219 (2000)


\bibitem[ShP]{ShP}
Shchepochkina I.,  Post, G., Explicit bracket in an exceptional
simple Lie superalgebra. Internat. J. Algebra Comput. 8 (1998),
no. 4, 479--495


\bibitem[Sh14]{Sh14}
Shchepochkina I., Five exceptional simple Lie superalgebras of
vector fields and their fourteen regradings.  Represent.  Theory,
v.  3, 1999, 3 (1999), 373--415; hep-th 9702121

\bibitem[ShE]{ShE}
Shchepochkina I., How to realize Lie algebras by vector fields.
Examples. ??

\bibitem[Sh2]{Sh2}
Shchepochkina I., How to realize Lie algebras by vector fields.
More examples. ??

\bibitem[St]{St}
Sternberg S., {\em Lectures on differential geometry}, Chelsey,
2nd edition, 1985

\bibitem[Y]{Y}
Yamaguchi K., Differential systems associated with simple graded
Lie algebras. Progress in differential geometry, Adv. Stud. Pure
Math., 22, Math. Soc. Japan, Tokyo, 1993, 413--494


\end{thebibliography}
\end{document}